\RequirePackage{ifpdf}
\ifpdf % We are running pdfTeX in pdf mode
\documentclass[pdftex]{sigma}
\else
\documentclass{sigma}
\fi

\usepackage{mathrsfs}
\usepackage{stmaryrd}
\usepackage{euscript}
\newcommand{\set}[1]{\left\{#1\right\}}

\usepackage{amssymb}

\DeclareMathOperator{\spann}{\dsr-span}
\DeclareMathOperator{\spannq}{\dsq-span}
\DeclareMathOperator{\zspann}{\dsz    -span}

\DeclareMathOperator{\ad}{ad}
\DeclareMathOperator{\bbase}{\mathscr{B}}

 \DeclareMathOperator{\Aut}{Aut}
\DeclareMathOperator{\derived}{\EuScript{D}(\mathfrak g)}
\DeclareMathOperator{\derivedg}{\EuScript{D}(G)}

\DeclareMathOperator{\ens}{\EuScript{A}}
\DeclareMathOperator{\enss}{\EuScript{B}}
        \newfont{\Sans}{cmss10 scaled\magstep1}
\newcommand{\sans}[1]{\hbox{\Sans {#1}}}
\def\hb#1{\hbox{#1}}

\def\cent{\mbox{\usefont{T1}{bch}{bx}{n}
  \selectfont c}}
  \def\centg{\mbox{\usefont{T1}{bch}{bx}{n}
  \selectfont C}}

\def\ker#1{\hb{ker}(#1)}

\def\hb #1{\hbox{#1}}

\def\hb#1{\hbox{#1}}

\def\ker#1{\hb{ker}(#1)}

\def\dim#1{\hb{dim}(#1)}

\def\lef({\left(}
\def\rig){\right)}

\def\dsr{{\mathbb R}}

\def\dsn{{\mathbb N}}
\def\dsq{{\mathbb Q}}
\def\dsz{{\mathbb Z}}
\let\dis=\displaystyle

\begin{document}

\allowdisplaybreaks

\renewcommand{\thefootnote}{$\star$}

\renewcommand{\PaperNumber}{020}

\FirstPageHeading

\ShortArticleName{Discrete Cocompact Subgroups of Nilpotent Lie
Groups}

\ArticleName{Discrete Cocompact Subgroups\\ of  the Five-Dimensional
Connected\\ and Simply Connected Nilpotent Lie Groups}

\Author{Amira GHORBEL and Hatem HAMROUNI}

\AuthorNameForHeading{A. Ghorbel and H. Hamrouni}

\Address{Department of Mathematics, Faculty of Sciences at Sfax,\\ Route
Soukra, B.P.~1171, 3000 Sfax, Tunisia}
\Email{\href{mailto:Amira.Ghorbel@fss.rnu.tn}{Amira.Ghorbel@fss.rnu.tn}, \href{mailto:hatemhhamrouni@voila.fr}{hatemhhamrouni@voila.fr}}

\ArticleDates{Received July 16, 2008, in f\/inal form February 09,
2009; Published online February 17, 2009}

\Abstract{The discrete cocompact subgroups of the f\/ive-dimensional
connected, simply connected
 nilpotent  Lie groups  are determined up to isomorphism. Moreover,
 we prove if $G=N\times A$ is a connected, simply connected,
 nilpotent Lie group with an Abelian factor $A$, then every uniform
 subgroup of $G$ is the direct product of a uniform subgroup of $N$
 and $\dsz^r$ where $r=\dim{A}$.}


\Classification{22E40}

\vspace{-2mm}

\section{Introduction}
To solve many problems in various branches of mathematics, for
example, in dif\/ferential geo\-metry (invariant geometric structures on
Lie groups), in spectral geometry, in physics (multi-dimensional
models of space-times) etc.\ (see \cite{Gorba1}), we are led to study
the explicit description of uniform subgroups (i.e., discrete
cocompact) of solvable Lie groups. Towards that purpose, we focus
attention to the description of uniform subgroups of a connected,
simply connected non Abelian nilpotent Lie groups satisfying the
rationality criterion of Malcev (see \cite{Malcev1}). At present we
have no progress  on this dif\/f\/icult problem although the
classif\/ication for certain groups is already exists. To attack this
problem, we will be in a f\/irst step interested  in giving an
explicit description of uniform subgroups of
 nilpotent Lie groups of dimension less or equal to $6$. \mbox{Another}
motivation  comes from \cite[p.~339]{Milnes100}. For the case
of  $3$ or $4$ dimensions, the classif\/ication is given
 in \cite{Auslander1,Pesce1,Milnes1}. For  $5$-dimensional nilpotent
 Lie groups,
  there
are eight real connected and simply connected non Abelian nilpotent
Lie groups which are $G_3\times \dsr^2$, $G_4\times \dsr$ and
$G_{5,i}$ for $1\leq i\leq 6$~\cite{Dixmier5}. The uniform
subgroups of $G_3\times \dsr^2$, $G_{5,1}$, $G_{5,3}$ and $G_{5,5}$
are determined respectively in \cite{Pesce1,Carolyn1,
Milnes2,Hamrouni1}. The aim of this paper is to
 complete the classif\/ication of the discrete cocompact subgroups of the Lie groups
$G_4\times \dsr$, $G_{5,2}$, $G_{5,4}$ and $G_{5,6}$.

This paper is organized as follows. In Section~\ref{section2}, we f\/ix some
notation which will be of use later and we record   few standard
facts about rational structures and uniform subgroups of a~connected, simply connected nilpotent Lie groups. Section~\ref{section3} is
devoted to the study of the rationality of certain subalgebras of a
given nilpotent Lie algebra with a rational structure. In Section~\ref{section4}, in order to determine the uniform subgroups of $G_4\times
\dsr$, we study a more general class when $G$ is a connected, simply
connected, nilpotent Lie group with an Abelian factor, that is
$G=N\times A$ where $A$ is an Abelian normal subgroup of $G$.
Theorem \ref{non-cancelable} shows that every uniform subgroup
$\Gamma$ of $G$ is a direct product of a uniform subgroup of $N$ and
$\dsz^r$ where $r=\dim{A}$. As an immediate application of this
result, we determine the uniform subgroups of $G_4\times \dsr$. We
determine afterwards  all uniform subgroups of the Lie groups
$G_{5,2}$, $G_{5,4}$ and $G_{5,6}$.

\section{Notations and basic facts}\label{section2}

The aim of this section is to give a brief review of certain results
from rational structures and uniform subgroups of connected and
simply connected nilpotent
 Lie groups which will be needed
later. The reader who is interested in detailed proof is referred to
standard texts \cite{Cor1, Raghunathan,Malcev1}.

\subsection{Rational structures and uniform subgroups}\label{section2.1}

   Let $G$ be a nilpotent, connected
and simply connected real Lie group and let $\mathfrak g$ be its Lie
algebra. We say that $\mathfrak g$ (or $G$) has a \textit{rational
structure} if there is a Lie algebra $\mathfrak g_\mathbb{Q}$ over
$\mathbb{Q}$ such that $\mathfrak g \cong \mathfrak g_\mathbb{Q} \otimes
\mathbb{R}$. It is clear that $\mathfrak g$ has a rational structure if
and only if $\mathfrak g$ has an $ \mathbb{R}$-basis
$\{X_1,\dotsc,X_n\}$ with rational structure constants.

 Let $\mathfrak g$ have a f\/ixed
rational structure given by $\mathfrak g_\mathbb{Q}$ and let $\mathfrak h$
be an $\mathbb{R}$-subspace of $\mathfrak g$. Def\/ine $\mathfrak
h_\mathbb{Q}= \mathfrak h \cap \mathfrak g_\mathbb{Q}$. We say that $\mathfrak
h$ is \textit{rational} if $\mathfrak h = \spann\set{\mathfrak
h_\mathbb{Q}}$, and that a connected, closed subgroup~$H$ of $G$ is
\textit{rational} if its Lie algebra $\mathfrak h$ is rational. The
elements of $\mathfrak g_\mathbb{Q}$ (or $G_\mathbb{Q}= \exp(\mathfrak
g_\mathbb{Q})$) are called \textit{rational elements} (or
\textit{rational points}) of $\mathfrak g$ (or $G$).

 A
discrete subgroup $\Gamma$ is called \textit{uniform} in $G$ if the
quotient space $G/\Gamma$ is compact. The homogeneous space
$G/\Gamma$ is called a \textit{compact nilmanifold}. A proof of the
next result can be found in Theorem~7 of \cite{Malcev1} or in
Theorem~2.12 of \cite{Raghunathan}.

\begin{theorem}[The Malcev rationality criterion]
Let $G$ be a simply connected nilpotent Lie group, and let $\mathfrak g$
be its Lie algebra. Then $G$ admits a uniform subgroup $\Gamma$ if
and only if $\mathfrak g$ admits a basis $\{X_1,\dotsc,X_n\}$ such that
\[
[X_i, X_j] = \dis\sum_{\alpha=1}^n c_{ij\alpha} X_\alpha\qquad
\mbox{for all} \quad i, \ j,
\] where the constants $c_{ij\alpha}$ are all
rational. (The $c_{ij\alpha}$ are called the structure constants of
$\mathfrak g$ relative to the basis $\set{X_1, \ldots,
X_n}$.)\end{theorem}

More precisely, we have, if $G$ has a uniform
subgroup $\Gamma$, then $\mathfrak g$ (hence $G$) has a rational
structure such that $\mathfrak g_\mathbb{Q}= \spannq\set{\log(\Gamma)}$.
Conversely, if $\mathfrak g$ has a rational structure given by some
$\mathbb{Q}$-algebra $\mathfrak g_\mathbb{Q}\subset \mathfrak g$, then $G$
has a uniform subgroup $\Gamma$ such that $\log(\Gamma) \subset
\mathfrak g_\mathbb{Q}$ (see~\cite{Cor1,Malcev1}).  If we
endow $G$ with the rational structure induced by a uniform subgroup
$\Gamma$ and if $H$ is a Lie subgroup of $G$, then $H$ is rational
if and only if $H\cap \Gamma$ is a uniform subgroup of $H$. Note
that the notion of rational depends on $\Gamma$.

\subsubsection{Weak and strong Malcev basis}

 Let $\mathfrak g$ be a
nilpotent Lie algebra and let $\bbase=\set{X_1, \ldots, X_n}$ be a
basis of $\mathfrak g$.
 We say that $\bbase$ is a weak (resp.\ strong) Malcev
basis for $\mathfrak g$ if  $\mathfrak g_i=\spann\set{X_1, \ldots, X_i}$ is
a subalgebras (resp.\ an ideal) of $\mathfrak g$ for each $1\leq i\leq n$
(see~\cite{Cor1}).

 Let $ \Gamma$ be a uniform
subgroup of $G$. A strong Malcev (or Jordan--H\"{o}lder) basis
$\{X_1,\dotsc,X_n\}$ for $\mathfrak g$ is said to be \textit{strongly
based on} $\Gamma$ if
\begin{equation*}%\label{strongly basis}
\Gamma= \exp (\mathbb{Z}X_1)\cdots  \exp (\mathbb{Z}X_n).
\end{equation*}
 Such a basis always exists (see \cite{Cor1,Matsushima1}).

 The \textit{lower central series}
(or the \textit{descending central series}) of $\mathfrak g$ is the
decreasing sequence of characteristic ideals of $\mathfrak g$ def\/ined
inductively as follows
\[
 \EuScript{C}^1(\mathfrak g) = \mathfrak g;\qquad \EuScript{C}^{p+1}
(\mathfrak g) = [\mathfrak g, \EuScript{C}^p(\mathfrak g)]\qquad (p \ge 1).
\] The
characteristic ideal $ \EuScript{C}^2(\mathfrak g) = [\mathfrak g, \mathfrak g]$
is called the derived ideal of the Lie algebra $\mathfrak g$ and denoted
by $ \derived$. Let $\derivedg = \exp(\derived)$, observe that we
have $\derivedg = [G, G]$.

 The Lie algebra $\mathfrak g$
is called $k$-\textit{step nilpotent Lie algebra} or
\textit{nilpotent Lie algebra of class} $k$ if there is an integer
$k$ such that
\[ \EuScript{C}^{k+1}(\mathfrak g) = \set {0}\qquad \mbox{and}\qquad
\EuScript{C}^k(\mathfrak g)\ne \set {0}.
\] We denote the center of $G$
by $Z(G)$ and the center of $\mathfrak g$ by $\mathfrak{z}(\mathfrak{g})$.

\begin{proposition}[\cite{Matsushima1,Cor1}]\label{rational-descending-central-series}  If $\mathfrak g$
has rational structure, all the algebras in the descending central
series are rational.
\end{proposition}
  Let $G$ be a group. The center $Z(G)$ of $G$ is a normal
subgroup. Let $\mathscr{C}_2(G)$ be the inverse image of $Z(G/Z(G))$
under the canonical projection $G \longrightarrow G/Z(G)$. Then
$\mathscr{C}_2(G)$ is normal in~$G$ and contains $Z(G)$. Continue
this process by def\/ining inductively: $\mathscr{C}_1(G)=Z(G) $ and
$\mathscr{C}_i(G)$ is the inverse image of
$Z(G/\mathscr{C}_{i-1}(G))$ under the canonical projection $G
\longrightarrow G/\mathscr{C}_{i-1}(G)$. Thus we obtain a sequence
of normal subgroups of $G$, called the ascending central series of
$G$.
\begin{proposition}[\cite{Cor1}]\label{rational-ascending-central-series} If $G$ is a
nilpotent Lie group with rational structure, all the algebras in the
ascending central series are rational. In particular, the center
$\mathfrak z (\mathfrak g)$ of $\mathfrak g$ is rational.
\end{proposition}

A proof of the next result can be found in Proposition 5.3.2 of
\cite{Cor1}.
\begin{proposition}
Let $\Gamma$ be uniform subgroup in a nilpotent Lie group $G$, and
let $H_1\subsetneqq H_2\subsetneqq\cdots\subsetneqq H_k =G$ be
rational Lie  subgroups of $G$. Let $\mathfrak h_1, \ldots, \mathfrak
h_{k-1}, \mathfrak h_k =\mathfrak g$ be the corresponding Lie algebras. Then
there exists a weak Malcev basis $\set{X_1,\ldots, X_n}$ for $\mathfrak
g$ strongly based on $\Gamma$ and passing through $\mathfrak h_1,
\ldots, \mathfrak h_{k-1}$. If the $H_j$ are all normal, the basis can
be chosen to be a strong Malcev basis.
\end{proposition}
  A rational structure on $\mathfrak g$ induces a rational structure
on the dual space $\mathfrak g^*$ (for further details, see \cite[Chaper~5]{Cor1}). If $\mathfrak g$ has a rational structure given by the
uniform subgroup $\Gamma$, a real linear functional $f \in \mathfrak
g^*$ is rational ($f\in \mathfrak g_\mathbb{Q}^*,  \mathfrak
g_\mathbb{Q}=\spannq\set{\log(\Gamma)}$) if $\langle f, \mathfrak
g_\mathbb{Q}\rangle \subset \mathbb{Q}$, or equivalently $\langle f,
\log(\Gamma) \rangle \subset \mathbb{Q}$.   Let
$\mathrm{{Aut}}(G)$ (respectively $\mathrm{{Aut}}(\mathfrak g)$)
 denote the group of automorphism of $G$ (respectively  $\mathfrak g $). If
 $\varphi\in \mbox{Aut}(G)$, $\varphi_*$ will denote the
 derivative of $\varphi$ at identity. The mapping $\mathrm{{Aut}}(G)
 \longrightarrow\mathrm{{Aut}}(\mathfrak g),\ \varphi \longmapsto
 \varphi_*$ is a groups isomorphism (since $G$ is simply
 connected).

\begin{theorem}[\protect{\cite[Theorem~5]{Malcev1}}]\label{extended}
Let $G_1$  and $G_2$ be connected simply connected nilpotent Lie
groups and
 $\Gamma_1$, $\Gamma_2$  uniform subgroups of $G_1$  and $G_2$. Any abstract group
isomorphism $f$ between $\Gamma_1$ and $\Gamma_2$ extends uniquely
to an
 isomorphism $\overline{f}$ of $G_1$ on $G_2$; that is, the
 following diagram
$$\begin{CD}
\Gamma_1 @>f>>\Gamma_2\\
@V{i}VV @VV{i}V\\
G_1@>>\overline{f}> G_2
\end{CD}$$
is commutative, where $i$ is the inclusion mapping.
\end{theorem}
  We  conclude this review with two results.

\subsection{Smith normal form}\label{section2.2}

 A commutative ring $\mathbf{R}$ with
identity $1_\mathbf{R}\ne 0$ and no zero divisors is called an
\textit{integral domain}. A principal ideal ring which is an
integral domain is called a \textit{principal ideal domain}.
\begin{theorem}[elementary divisors theorem]\label{Smith-normal-form} If
$A$ is an $n\times n$ matrix of rank $r>0$ over a~principal ideal
domain $\mathbf{R}$, then $A$ is equivalent to a matrix of the form
\[\left(
    \begin{array}{cc}
      L_r & 0 \\
      0 & 0 \\
    \end{array}
  \right),
  \]
where $L_r$ is an $r\times r$ diagonal matrix with nonzero diagonal
entries $d_1, {\ldots}, d_r$ such that $d_1|d_2| {\cdots}|d_r.\!$
\end{theorem}
The notation $d_1\vert d_2\vert \cdots \vert d_r$ means $d_1$
divides $d_2$, $d_2$ divides $d_3$, etc. The elements $d_1, \ldots,
d_r$ are called the elementary divisors of $A$. See
\cite{Hungerford} for a more precise result.

\subsection{Hermite normal form}\label{section2.3}

\begin{definition}[Hermite normal form] A matrix $(a_{ij})\in
\mathrm{Mat}(m, n, \dsr)$ with $m\leq n$ is in Hermite normal form
if the following conditions are satisf\/ied:
\begin{itemize}\itemsep=0pt
  \item[$(1)$] $a_{ij}=0$  for   $  i>j$;
  \item[$(2)$] $a_{ii}>0$ for $i=1, \ldots, m$;
  \item[$(3)$] $0\leq a_{ij}< a_{ii}$ for $i<j$.
\end{itemize}

\end{definition}
\begin{theorem}[Hermite 1850]\label{Hermite} For every matrix $A \in \mathrm{Mat}(m, n, \dsr)$
with $\mathrm{rank}(A) =m\leq n$, there is a matrix $T\in
\mathrm{GL}(n, \dsz)$, so that $AT$ is in Hermite normal form. The
Hermite normal form $AT$ is unique.
\end{theorem}

 \section{Rationality of certain  subalgebras}\label{section3}

 The following generalizes the Proposition 5.2.4 of \cite{Cor1}.
 \begin{proposition} Let $G$ be a simply connected nilpotent Lie group and
 $\Gamma\subset G$ a uniform subgroup. Let $\Delta\subset \Gamma$
 be a finite set and let $\centg(\Delta)$ denote the centralizer of
 $\Delta$ in $G$. Then $\centg(\Delta)$ is rational.
 \end{proposition}
 \begin{proof} This follows immediately from Lemma 1.14, Theorem 2.1
 of \cite{Raghunathan} and Theorem 4.5 of~\cite{Onishchik1}.~~
 \end{proof}

  Next, we prove the following proposition which will play an
important role below.

 \begin{proposition}\label{rat-centralizer} Let $\Gamma$ be a uniform
 subgroup of a nilpotent Lie group $G=\exp(\mathfrak g)$.
Let $H = \exp(\mathfrak h)$ be a rational subgroup of  $G$. Then the
centralizer  $ \cent(\mathfrak h)$ of $\mathfrak h$ in $\mathfrak g$ is
rational.
\end{proposition}
\begin{proof} Let $\set{e_1, \ldots, e_n}$ be a weak Malcev basis for $\mathfrak
g$ strongly based on  $ \Gamma$ passing through $\mathfrak h$ (see
\cite[Proposition 5.3.2]{Cor1}). We note $\mathfrak h = \spann \set{e_1,
\ldots, e_p}$. Let $\set{e_1^*, \ldots, e_n^*}$ be the dual basis of
$\set{e_1, \ldots, e_n}$. For $i=1, \ldots, n$ and $j=1, \ldots, p$,
we def\/ine
\[
f_{ij}: \ \mathfrak g \longrightarrow \dsr;\qquad X \longmapsto
\langle e_i^*, [X, e_j]\rangle.
\] The functionals $f_{ij}$ are
rational. Then the kernels $\ker{f_{ij}}$ are rational subspaces in
$\mathfrak g$ (see \cite[Lemma 5.1.2]{Cor1}). On the other hand, it is
easy to see that
\[
\cent(\mathfrak h) = \bigcap_{{1\leq i\leq n }\atop
{1\leq j\leq p }} \ker{f_{ij}}.
\]
 Therefore, we conclude from Lemma~5.1.2 of \cite{Cor1} that $\cent(\mathfrak h)$ is rational.
 %This completes the proof.
\end{proof}

\section[Uniform subgroups of   nilpotent Lie
group with an Abelian factor.
Uniform subgroups of $G_4\times \dsr$]{Uniform subgroups of   nilpotent Lie
group\\ with an Abelian factor.
Uniform subgroups of $\boldsymbol{G_4\times \dsr}$} \label{section4}

The aim of this section is
 to describe the classif\/ication of the uniform subgroups of $G_4\times \dsr$.
Let us begin with a more general situation. First, we introduce the
following def\/inition.
\begin{definition}[\cite{Lauret1,Gorba1}]
 An \textit{Abelian
factor} of a  Lie algebra $\mathfrak g$ is an Abelian ideal $\mathfrak a$
for which there exists an ideal $\mathfrak n$ of $\mathfrak g$ such that
$\mathfrak g = \mathfrak n \oplus \mathfrak a$ (i.e., $[\mathfrak n , \mathfrak a] =
\set{0}$).
\end{definition}
Let $m(\mathfrak g)$ denote the maximum dimension over all Abelian
factors of $\mathfrak g$. If $\mathfrak z(\mathfrak g)$ is the center of $\mathfrak
g$ then the maximal Abelian factors are precisely the linear direct
complements of $\mathfrak z(\mathfrak g)\cap \derived$ in $\mathfrak z(\mathfrak
g)$, that is, those subspaces $\mathfrak a\subset \mathfrak z(\mathfrak g)$ such
that $\mathfrak z(\mathfrak g)= \mathfrak z(\mathfrak g)\cap \derived \oplus \mathfrak
a$. Therefore
\[
m(\mathfrak g) = \dim{\mathfrak z(\mathfrak g)} - \dim{\mathfrak
z(\mathfrak g)\cap \derived}.
\]

Let $\mathfrak g$ be a nilpotent Lie algebra and $\mathfrak a$ an Abelian
factor of $\mathfrak g$. Let $\mathfrak g = \mathfrak n\oplus \dsr^r$ be any
decomposition in ideals of $\mathfrak g$. The next simple lemma
establishes a relation between the derived algebra of $\mathfrak g$ and
the derived algebra of $\mathfrak n$. Its value will be immediately
apparent in the proof of Theorem~\ref{non-cancelable}.

\begin{lemma} With the above notation, we have
\begin{equation*}
\derived = \EuScript{D}(\mathfrak n).
\end{equation*}
\end{lemma}
\begin{proof} Let $X_1, X_2\in \mathfrak g$. Write $X_i=a_i+b_i$, $i=1,2$,
where $a_i\in \mathfrak a$ and $b_i\in \mathfrak n$. We calculate
\begin{gather*}
  [X_1, X_2]  =  [a_1+b_1, a_2+b_2]  = [b_1, b_2]
\end{gather*}
since $\mathfrak a\subset \mathfrak z(\mathfrak g)$. Therefore $\derived =
\EuScript{D}(\mathfrak n)$. The proof of the lemma is complete.
\end{proof}

In the sequel, the symbol $\simeq$ denotes
 abstract group isomorphism.
\begin{theorem}\label{non-cancelable}
 Let $\mathfrak g$ be a  nilpotent Lie algebra with maximal Abelian factor of dimension
 $m(\mathfrak
g)=r$ and let $\mathfrak g = \mathfrak n\oplus \dsr^r$ be any decomposition
in ideals, that is $\dsr^r$ is a maximal Abelian factor of $\mathfrak
g$. Then  we have the following:
\begin{itemize}\itemsep=0pt
  \item[$(1)$] The group $G=\exp(\mathfrak g)$ admits a uniform subgroup if and only if
$N=\exp(\mathfrak n)$ admits a~uniform subgroup.
  \item[$(2)$] If $\Gamma$ is a uniform subgroup of $G$, then there exists a
uniform subgroup $H$ of $N$ such that
\begin{equation*}
\Gamma \simeq H\times \dsz^r.
\end{equation*}
\end{itemize}
\end{theorem}
\begin{proof} Let $\mathfrak a = \dsr^r$.

 $(1)$ If $H$ is a uniform subgroup of $N$ then it is clear that
$H\times \dsz^r$ is a~uniform subgroup of~$G$. Conversely, we
suppose that $G$ admits a uniform subgroup $\Gamma$. Let $\set{e_1,
\ldots, e_n}$ be a strong Malcev basis for $\mathfrak g$ strongly based
on $\Gamma$ passing through $\derived$ and $\derived +\mathfrak z(\mathfrak
g)$. Put
\begin{equation*}
\derived = \spann\set{e_1, \ldots, e_q}
\end{equation*}
and
\begin{equation*}
\derived+ \mathfrak z(\mathfrak g) = \spann\set{e_1, \ldots, e_q, e_{q+1},
\ldots, e_{q+r}}.
\end{equation*}
For every $i=q+r+1, \ldots, n$, we note
\begin{equation*} e_i = u_i+v_i,
\end{equation*}
where $u_i\in \mathfrak n$ and $v_i\in \mathfrak a$. It is clear that
$\set{e_1, \ldots, e_q, u_{q+r+1}, \ldots, u_n}$ is a basis for
$\mathfrak n$ and has the same structure constants as $\set{e_1, \ldots,
e_n}$. By the Malcev rationality criterion, $N$ admits a~uniform
subgroup.

 $(2)$ Let $\Gamma$ be a uniform subgroup of $G$. Since $\mathfrak
z(\mathfrak g)$ and $\derived$ are rational then there exists  a~strong
Malcev basis $\set{e_1, \ldots, e_n}$ for $\mathfrak g$ strongly based
on $ \Gamma$ passing through $\mathfrak z(\mathfrak g)\cap\derived$, $\mathfrak
z(\mathfrak g)$ and $\derived+ \mathfrak z(\mathfrak g)$. Put
\[
\mathfrak
z(\mathfrak g)\cap\derived =\spann\set{e_1, \ldots, e_p}.
\]
 Since
$\dim{\mathfrak z(\mathfrak g)/\mathfrak z(\mathfrak g)\cap\derived} = m(\mathfrak g)
=r$, then
\[
\mathfrak z(\mathfrak g) =\spann\set{e_1, \ldots,
e_p, \ldots, e_{p+r}}.
\] Since $\mathfrak z(\mathfrak g)$ is Abelian, we can
assume that $(e_{1}, \ldots, e_p)$ are independent modulo $\mathfrak a$,
i.e., they span a complement of $\mathfrak a$ in $\mathfrak z(\mathfrak g)$. Let
for every $i=p+1, \ldots, n$
\[
e_i= a_i+b_i,
\]
where $a_i\in \mathfrak a$ and $b_i \in \mathfrak n$. On the other hand, we
have $\derived+ \mathfrak z(\mathfrak g)=\derived \oplus \mathfrak a$, then
$\mathfrak a = \spann\set{a_{p+1}, \ldots, a_{p+r}}$. Let $\Phi_*: \mathfrak
g \longrightarrow \mathfrak g$  be the function def\/ined by
\[
e_i \longmapsto \left\{
                    \begin{array}{ll}
                      e_i, & \hbox{if  } 1\leq i\leq p, \\
a_i, & \hbox{if  } p+1\leq i\leq p+r,\\
                      b_i, & \hbox{if  } p+r+1\leq i\leq n.
                    \end{array}
                  \right.
\]
We will show that $\Phi_*$ is a Lie automorphism of $\mathfrak g$. In
fact, if $i\leq p+r$ or $j\leq p+r$, we have $\Phi_*([e_i, e_j]) =
[\Phi_*(e_i), \Phi_*(e_j)] =0$. Next, we suppose that $i, j > p+r$.
It is easy to verify  that
\[
[\Phi_*(e_i), \Phi_*(e_j)] = [b_i, b_j].
\]
On the other hand, we have
\[
[e_i, e_j] = [a_i+b_i, a_j+b_j] = [b_i, b_j]
\]
and hence
\[
\Phi_*([e_i, e_j]) = [b_i, b_j].
\]
Consequently, we obtain
\[
\Phi_*([e_i, e_j]) = [\Phi_*(e_i), \Phi_*(e_j)].
\]
 Then
\begin{gather*}
  \Gamma \simeq \Phi(\Gamma)  =  \prod_{i=1}^p \exp(\dsz
e_i) \prod_{i=p+1}^{p+r} \exp(\dsz a_i) \prod_{i={p+r+1}}^n
\exp(\dsz
b_i) \\
 \phantom{\Gamma \simeq \Phi(\Gamma)}{}    =  \prod_{i=1}^p \exp(\dsz
e_i) \prod_{i={p+r+1}}^n \exp(\dsz b_i) \prod_{i=p+1}^{p+r}
\exp(\dsz a_i).
\end{gather*}
Evidently,  $H = \prod\limits_{i=1}^p \exp(\dsz e_i) \prod\limits_{i={p+r+1}}^n
\exp(\dsz b_i)$ is a uniform subgroup of $N$ and
$K=\prod\limits_{i=p+1}^{p+r} \exp(\dsz a_i)$ is a uniform subgroup of
$A=\exp(\mathfrak a)$. It follows that
\[
\Gamma
   \simeq H\times \dsz^r.
\]
 This completes the  proof of the theorem.
\end{proof}
  A consequence of  the above theorem, we deduce the uniform
subgroups of $G_4\times \dsr$. We recall that the Lie algebra $\mathfrak
g_4$ of $G_4$ is spanned by the vectors
  $X_1, \ldots, X_4$
   such that the only non vanishing
brackets are
\[
 [X_4, X_3] = X_2,\qquad [X_4, X_2] =X_1.
\]
\begin{corollary} Let $\set{e}$ be the canonical
basis of $\dsr$. Every uniform subgroup $\Gamma$ of $G_4\times \dsr$
has the following form
\[ \Gamma \simeq \exp(\dsz e)
\exp(\mathbb{Z}X_1) \exp(p_1\mathbb{Z}X_2) \exp\left(\mathbb{Z}\left(p_1p_2X_3
- \frac{p_3}{2}X_2\right)\right) \exp(\mathbb{Z}X_4),
\] where $p_1$, $p_2$, $p_3$ are
integers satisfying $p_1 >0$, $p_2>0$, $p_1p_2+p_3 \in 2\mathbb{Z}$ and
$0\leq p_3 < 2 p_1$. Furthermore, different choices for the $p$'s
give non isomorphic subgroups.
\end{corollary}
\begin{proof}  The proof
follows from Theorem \ref{non-cancelable} and \cite[Proposition 4.1]{Hamrouni1} (see also \cite[Theorem~1]{Milnes1}).\end{proof}

\begin{remark}
We note that Proposition B.1 of \cite{Pesce1} becomes a direct
consequence of Theorem \ref{non-cancelable} and Theorem 2.4 of~\cite{Carolyn1}.
\end{remark}

\section[Uniform subgroups of  $G_{5, 2}$]{Uniform subgroups of  $\boldsymbol{G_{5, 2}}$}\label{section5}

Let $\mathfrak g_{5,2}$ be the two-step nilpotent Lie algebra with basis
\[
\bbase=\set{X_1, \ldots, X_5}
\]
and non-trivial Lie  brackets def\/ined by
\begin{equation}\label{crochet5-2}
[X_5, X_4] = X_2,\qquad  [X_5, X_3] =X_1.\end{equation}
Let $G_{5,2}$ be the corresponding connected and simply connected
nilpotent Lie group. First, we introduce the following set:

\begin{equation*}\label{D2}
 \mathscr{D}_2 = \{ r=(r_1, r_2)\in
(\mathbb{N}^*)^2:\ r_1 \mbox{ divides } r_2\}.
\end{equation*}
\begin{theorem} \qquad {}
\begin{itemize}\itemsep=0pt
\item[$1.$]  Let $r=(r_1, r_2)\in \mathscr{D}_2$. Then
\begin{equation*}\label{uniform-g5-2}
\Gamma_r = \exp\left( \frac{1}{r_1}\dsz X_1\right) \exp\left(\frac{1}{r_2} \dsz X_2\right)
\exp(\dsz X_3) \exp(\dsz X_4) \exp(\dsz X_5)
\end{equation*}
 is a  uniform
subgroup of $G_{5,2}$.
 \item[$2.$] If $\Gamma$ is a uniform subgroup of $G_{5,2}$, then there exist
 $r\in \mathscr{D}_2$
 and $\Phi \in \mathrm{{Aut}}(G_{5,2})$ such that
  $\Phi(\Gamma)= \Gamma_r$.
\item[$3.$] For $r$ and
$s$  in $\mathscr{D}_2$, $\Gamma_r$ and $\Gamma_s$ are isomorphic
groups if and only if $r=s$. \end{itemize}
\end{theorem}

\begin{proof} The proof of this theorem will be achieved through a
sequence of partial results.

\begin{lemma} The group $ \mathrm{Aut}(\mathfrak g_{5,2})$ is the set of
all matrices of the form
\[
\left(
    \begin{array}{ccc}
      \alpha A & B & u \\
      0 & A & v \\
      0 & 0 & \alpha \\
    \end{array}
  \right),
  \]
where $A\in \mathrm{GL}(2, \dsr), B\in \mathrm{Mat}(2, \dsr),
\alpha\in \dsr^*$ and $u, v\in \dsr^2$.
\end{lemma}
\begin{proof} It is easy to verify that $\spann\set{X_1, \ldots,
X_4}$ is the unique one co-dimensional Abelian ideal of $\mathfrak g$.
Therefore any automorphism $A$ of $\mathfrak g$ leaves invariant the
subalgebras $\spann\set{X_1, X_2}$ and $\spann\set{X_1, \ldots,
X_4}$. The remainder of the proof follows from (\ref{crochet5-2}).
\end{proof}

\begin{lemma}\label{main-result}
Every uniform subgroup $ \Gamma$ of $G_{5, 2}$ has the following
form: there are   integers $p$, $q$ and $\alpha  $ satisfying $p, q
>0$ and $0\leq \alpha < q $ such that
\begin{equation*}\label{uniform-G(5,2)}
\Gamma \simeq \Gamma(p, q, \alpha)= \exp\left( \frac{1}{p}\dsz X_1\right) \exp\left(
\dsz\left(\frac{1}{q} X_2-\frac{\alpha}{pq}X_1\right)\right) \exp(\dsz X_3) \exp(\dsz
X_4) \exp(\dsz X_5).
\end{equation*}
\end{lemma}
\begin{proof}
Let $G=G_{5,2}$, $\mathfrak g = \mathfrak g_{5,2}$ and let $\Gamma$ be a
uniform subgroup of $G$.
 Let
\[
\Omega = \left\{l=\sum_{i=1}^{5}l_i X_i^*\in \mathfrak g^*:\ l_1 \ne
0\right\}
\] be the layer of the generic coadjoint orbits (see
\cite[Chapter 3]{Cor1}). As $\mathfrak g_\mathbb{Q}^* $  is dense in
$\mathfrak g^*$ and $\Omega$ is a non-empty Zariski open set in $\mathfrak
g^*$ then $\mathfrak g_\mathbb{Q}^*\cap \Omega \ne \varnothing$. Let $
l\in \mathfrak g_\mathbb{Q}^*\cap \Omega$. Since $l $ is rational then
by Proposition 5.2.6 of \cite{Cor1}, the radical $\mathfrak g(l)$  of
the skew-symmetric bilinear form $B_l$ on $\mathfrak g$ def\/ined by $
B_l(X, Y) = \langle l, [X, Y]\rangle\ (X, Y \in \mathfrak g)$, is also
rational. It follows by Proposition \ref{rat-centralizer} that
$\cent(\mathfrak g(l))$ is also rational subalgebra of $\mathfrak g$. As
\[
\mathfrak g(l) = \spann\set{X_1, X_2, l_2 X_3-l_1 X_4}
\] then a simple
calculation shows that
\[
\cent(\mathfrak g(l)) = \spann\set{ X_1,
\ldots, X_4}.
\] On the other hand, the derived ideal $[\mathfrak g,
\mathfrak g] = \spann\set{X_1, X_2}$ is also rational subalgebra in
$\mathfrak g$ (see \cite[Corollary 5.2.2]{Cor1}). It follows by
Proposition 5.3.2 of \cite{Cor1} that there exists   a strong Malcev
basis $\set{Y_1, \ldots, Y_5}$  of $\mathfrak g$ strongly based on
$\Gamma$ passing through $[\mathfrak g, \mathfrak g]$ and $ \spann\set{ X_1,
\ldots, X_4}$. Then, we have $ Y_1 = a_{11} X_1 + a_{12} X_2$,
    $Y_2 = a_{21} X_1 + a_{22} X_2$,
    $Y_i = \sum\limits_{j=1}^4 a_{ij} X_j$  $(i=3,4)$ and
    $Y_5= \sum\limits_{j=1}^5 a_{5j} X_j$,
where $a_{ij} \in \mathbb{R}$. The mapping $(\Phi_1)_*: \mathfrak g
\longrightarrow\mathfrak g$ def\/ined by
\begin{gather*}
(\Phi_1)_*(Y_i)  =  Y_i  \qquad (i=1,2), \\
  (\Phi_1)_*(Y_3)  =  \frac{1}{a_{55}}(a_{33} X_3+ a_{34} X_4), \\
  (\Phi_1)_*(Y_4)  =  \frac{1}{a_{55}}(a_{43} X_3+ a_{44} X_4), \\
  (\Phi_1)_*(Y_5)  =  X_5
\end{gather*}
 is a Lie algebra automorphism. Then
\begin{equation*}\Gamma \simeq \Gamma_1 = \exp(\dsz Y_1)  \exp(\dsz
Y_2)  \exp((\Phi_1)_*(Y_3))\exp((\Phi_1)_*(Y_4)) \exp(\dsz
X_5).\end{equation*} On the other hand, let the  Lie algebra
automorphism $(\Phi_2)_*: \mathfrak g \longrightarrow\mathfrak g$ def\/ined by
\begin{gather*}
(\Phi_2)_*(X_5) = X_5,\qquad (\Phi_2)_*(X_4)  = \frac{a_{43}}{a_{55}}
X_3+ \frac{a_{44}}{a_{55}} X_4,\qquad  (\Phi_2)_*(X_3)  =
\frac{a_{33}}{a_{55}} X_3+ \frac{a_{34}}{a_{55}} X_4.
\end{gather*} It follows that there exist $a, b, c, d\in \dsr$
such that
\begin{equation*}\Gamma_1 \simeq \Phi_2^{-1}(\Gamma_1) = \Gamma_2
 = \exp(\dsz(aX_1+bX_2))  \exp(\dsz(cX_1+dX_2))
  \prod_{i=3}^5\exp(\dsz X_i).\end{equation*} Next,
since $\exp(X_5) \exp(X_{3}) \exp(-X_5)\in \Gamma_2$, then
$\exp(X_1)\in \Gamma_2$ and hence the ideal $\spann\!\set{X_1}$ is
rational relative to $\Gamma_2$. It follows that there exist $x,
y_1, y_2\in \dsr$ such that
\begin{equation*}\Gamma_2
 = \exp(\dsz(xX_1))  \exp(\dsz(y_1X_1+y_2X_2))
  \prod_{i=3}^5\exp(\dsz X_i).\end{equation*}
Also, we use that $\exp(X_1)$ belongs to $\Gamma_2$, we deduce that
there exists $p\in \dsn^*$ such that $x= \frac{1}{p}$. Similarly,
since $\exp(X_5) \exp(X_{4}) \exp(-X_5)\in \Gamma_2$, then
$\exp(X_2)\in \Gamma_2$. Therefore there exists $(q, m)\in
\dsn^*\times\dsz$ such that $y_2=\frac{1}{q}$ and
$y_1=-\frac{m}{pq}$. Then
\begin{equation*} \Gamma_2
 = \exp\left(\frac{1}{p}\dsz X_1\right)  \exp\left(\dsz\left(\frac{1}{q}X_2-
\frac{m}{pq} X_1\right)\right)
  \prod_{i=3}^5\exp(\dsz X_i).\end{equation*} Finally, we observe that, if $\alpha$ is the
remainder of the division of $m$ by $q$, then
\begin{gather*} \Gamma_2
 = \exp\left(\frac{1}{p}\dsz X_1\right)  \exp\left(\dsz\left(\frac{1}{q}X_2-
\frac{\alpha}{pq} X_1\right)\right)\prod_{i=3}^5\exp(\dsz X_i).\tag*{\qed}
\end{gather*}
\renewcommand{\qed}{}
\end{proof}
\begin{lemma}\label{etapo2}
A necessary and sufficient condition that two subgroups  $\Gamma(p,
q, \alpha)$ and $ \Gamma(p', q', \alpha')$ are isomorphic is that
there exist  $A, B\in \mbox{GL}(2, \dsz)$ such that
\begin{equation}\label{eq:12} \left(
      \begin{array}{cc}
        p' & \ \ \alpha'\\
        0 & \ \ q' \\
      \end{array}
    \right) = A \left(
      \begin{array}{cc}
        p & \ \ \alpha\\
        0 & \ \ q \\
      \end{array}
    \right) B.\end{equation}
\end{lemma}
\begin{proof}
Let $\Gamma(p, q, \alpha) \simeq \Gamma(p', q', \alpha')$. It is
well known that any abstract isomorphism of $\Gamma(p, q, \alpha)$
onto $\Gamma(p', q', \alpha')$ is the restriction of an automorphism
$\Phi$ of $G$ (see \cite[Theorem 5, p.~292]{Malcev1}). Since
$\spann\set{X_1, \ldots, X_4}$ is the unique one co-dimensional
Abelian ideal of $\mathfrak g$, then we can  suppose  that $\Phi_*(X_5)
= X_5$, $\Phi_*(X_3)= a X_3+ b X_4$ and $\Phi_*(X_4)= c X_3+ d X_4$
($a, b, c, d\in \dsr$). We deduce that
\[
\zspann\set{a X_3+ b X_4, c
X_3+ d X_4}= \zspann\set{X_3, X_4}
\] and
\begin{gather*}
\zspann\set{\frac{1}{p}(a X_1+ b X_2), \frac{1}{q}(c X_1+ d
X_2)-\frac{\alpha}{pq}(a X_1+ b X_2)}\\
\qquad{}=\zspann\set{\frac{1}{p'}
X_1, \frac{1}{q'} X_2-\frac{\alpha'}{p'q'} X_1}.
\end{gather*}
 Consequently, we
obtain
\[
\left(
    \begin{array}{cc}
      a & c \\
      b & d \\
    \end{array}
  \right)\in \mbox{GL}(2, \dsz) \qquad \mbox{and}\qquad
  \left(
      \begin{array}{cc}
        p' & \alpha' \\
        0 & q' \\
      \end{array}
    \right)
    \left(
      \begin{array}{cc}
        a & c \\
        b & d \\
      \end{array}
    \right)
    \left(
      \begin{array}{cc}
        \frac{1}{p} & \frac{-\alpha}{pq}\vspace{1mm}\\
        0 & \frac{1}{q} \\
      \end{array}
    \right) \in \mbox{GL}(2, \dsz).
  \]
    Conversely, suppose that there exist $A, B\in \mbox{GL}(2,
    \dsz)$ satisfy the relation (\ref{eq:12}).
     The mapping $\phi_*$ def\/ined by
\[
\mbox{Mat}(\phi_*, \bbase) = \left(
      \begin{array}{ccc}
        B^{-1} & 0 &0\\
        0 & B^{-1} &0 \\
        0&0&1\\
      \end{array}
    \right).
  \]
     belongs to $\Aut(\mathfrak g)$ and  it is clear that $\phi(\Gamma(p, q, \alpha))
     = \Gamma(p', q', \alpha')$. The lemma is  completely proved.
     \end{proof}
Finally, an appeal to Lemma \ref{etapo2} and Smith normal form
completes the proof of Lemma~\ref{main-result}.
\end{proof}

\section[Uniform subgroups of  $G_{5, 4}$]{Uniform subgroups of  $\boldsymbol{G_{5, 4}}$}\label{section6}

Let $\mathfrak g_{5,4}$ be the three-step nilpotent Lie algebra with
basis
\[
\bbase=(X_1, \ldots, X_5)
\]
 with Lie brackets are given by
\begin{equation}\label{crochet4}
[X_5, X_4]=X_3,\qquad [X_5, X_3]=X_2,\qquad [X_4, X_3]=X_1.
\end{equation}
and the non-def\/ined brackets being equal to zero or obtained by
antisymmetry. Let $G_{5,4}$ be the corresponding connected and
simply connected nilpotent Lie group. For $u\in \dsz$, let
\[
\mathrm{\mathbf{A}}(u)= \left(
                            \begin{array}{cccc}
                              1 & 0 & \frac{u}{2} & 0 \\
                              0 & 1 & 0 & 0 \\
                              0 & 0 & 1 & 0 \\
                              0 & 0 & 0 & 1 \\
                            \end{array}
                          \right)
\left(\begin{array}{cccc}
1 & 0&0&0\\
0 & 1&1&\frac{1}{2} \\
0&0&1&1\\
0&0&0&1
\end{array} \right)= \left(\begin{array}{cccc}
1 & 0&\frac{u}{2}&\frac{u}{2}\\
0 & 1&1&\frac{1}{2} \\
0&0&1&1\\
0&0&0&1
\end{array} \right)
\]
and let
\[
\mathrm{\mathbf{B}} = \left(\begin{array}{ccc}
1 & 0&1\\
0 & 1&0\\
0&0&1\\
\end{array} \right).
\]

We denote by $\ens_4$  the subset of $\mbox{Mat}(4,
\mathbb{Z})$ consisting of all matrices of the form
\begin{equation*}%\label{matrice forme}
\llbracket D, m  \rrbracket = \begin{pmatrix}
  D & 0 \\
  0 & m
\end{pmatrix}
\end{equation*}
satisfying
\begin{equation}\label{Mariem1}
\llbracket D, m  \rrbracket^{-1} \mathrm{\mathbf{A}}(m) \llbracket
D, m \rrbracket \in \mbox{SL}(4, \mathbb{Z})
\end{equation}
and
\begin{equation}\label{Mariem2}
D^{-1} \mathrm{\mathbf{B}} D \in \mbox{SL}(3, \mathbb{Z}),
\end{equation}
where $m \in  \mathbb{N}^*$, the block matrix $D =(\alpha_{ij};
1\leq i, j\leq 3)$ is an upper-triangular integer invertible matrix
and $\mbox{SL}(3, \mathbb{Z})$ is the set of all integer matrices
with determinant $1$.

\begin{proposition} If $\llbracket D, m  \rrbracket\in \ens_4$
and if $H$ is the Hermite normal form of $D$, then $\llbracket H, m
\rrbracket\in \ens_4$.
\end{proposition}
\begin{proof} Let $H$ be the Hermite normal form of $D$ and let $T\in \mathrm{GL}(3, \dsz)$
such that $H = D T$. It is clear that $\llbracket H, m \rrbracket$
is the Hermite normal form of $\llbracket D, m  \rrbracket$ and
\[
\llbracket H, m
\rrbracket = \llbracket D, m \rrbracket \left(
                                          \begin{array}{cc}
                                            T & 0 \\
                                            0 & 1 \\
                                          \end{array}
                                        \right).
\]
Consequently \begin{equation*} \llbracket H, m \rrbracket ^{-1}
\mathrm{\mathbf{A}}(m) \llbracket H, m \rrbracket = \left(
                                          \begin{array}{cc}
                                            T & 0 \\
                                            0 & 1 \\
                                          \end{array}
                                        \right)^{-1} \llbracket D, m \rrbracket ^{-1}
\mathrm{\mathbf{A}}(m) \llbracket D, m \rrbracket \left(
                                          \begin{array}{cc}
                                            T & 0 \\
                                            0 & 1 \\
                                          \end{array}
                                        \right).
\end{equation*}
As $\llbracket D, m \rrbracket ^{-1} \mathrm{\mathbf{A}}(m)
\llbracket D, m \rrbracket\in \mbox{SL}(3, \mathbb{Z})$, then
$\llbracket H, m \rrbracket ^{-1} \mathrm{\mathbf{A}}(m) \llbracket
H, m \rrbracket \in \mbox{SL}(3, \mathbb{Z})$. The second condition
(\ref{Mariem2}) is shown similarly.
\end{proof}
Next, we def\/ine
\[
\enss_4= \Big\{\llbracket D, m  \rrbracket\in
\ens_4:\ D \mbox{ is in  Hermite normal form}\Big\}.
\]   We are
ready to formulate our result.
\begin{theorem}\label{uniform(5,4)} With the above notations, we have
\begin{itemize}\itemsep=0pt
\item[$1.$]  If $\llbracket D, m  \rrbracket\in \enss_4$, then
\begin{equation*}\label{uniform-g5-2}
\Gamma_{\llbracket D, m  \rrbracket} = \exp( \dsz \epsilon_1)
\exp(\dsz \epsilon_2) \exp(\dsz \epsilon_3) \exp(\dsz \epsilon_4)
\exp(\dsz \epsilon_5),
\end{equation*}
 where the vectors $\epsilon_j\ (1\leq j\leq 4)$ are the column vectors of
 $\llbracket D, m  \rrbracket$ in the basis $(X_1, \ldots, X_4)$ and
 $\epsilon_5=X_5$, is a discrete uniform
subgroup of $G_{5,4}$.
 \item[$2.$] If $\Gamma$ is a uniform subgroup of $G_{5,4}$, then there exist
 $\llbracket D, m  \rrbracket\in \enss_4$ and $\Phi \in \mathrm{{Aut}}(G_{5,4})$ such that
  $\Phi(\Gamma)= \Gamma_{\llbracket D, m  \rrbracket}$.
 \end{itemize}
\end{theorem}

 In the proof of the theorem  we need the following lemma.
\begin{lemma} The group $\mathrm{{Aut}}(\mathfrak g_{5,4})$ is the set of
all invertible matrices of the form
\begin{equation*}
\left(
                                           \begin{array}{ccccc}
                                             a_{44}\delta & a_{45}\delta & a_{45}a_{34}
-a_{35}a_{44} & a_{14} & a_{15} \\
                                             a_{54}\delta & a_{55}\delta & a_{55}a_{34}
-a_{35}a_{54} & a_{24} & a_{25} \\
                                             0 & 0 & \delta & a_{34} & a_{35} \\
                                             0 & 0 & 0 & a_{44} & a_{45}\\
                                             0 & 0 & 0 & a_{54} & a_{55} \\
                                           \end{array}
                                         \right),
\end{equation*}
where $\delta= a_{55}a_{44} -a_{45}a_{54}$.
\end{lemma}
\begin{proof}
Any automorphism $A$ of $\mathfrak g_{5,4}$ leaves invariant the
subalgebras $\spann\set{X_1, X_2}$ and $\spann\set{X_1, X_2, X_3}$.
The remainder of the proof follows from (\ref{crochet4}).
\end{proof}

\begin{proof}[Proof of Theorem \ref{uniform(5,4)}] Let $G=G_{5,4}$
and $\mathfrak g = \mathfrak g_{5,4}$.
 Assertion $1$ is obvious. To prove the second, let $\Gamma$ be a uniform subgroup
of $G$. Since  the ideals
\begin{equation*}\mathfrak z(\mathfrak g) =\spann\set{X_1, X_2} \qquad \mbox{and} \qquad \derived
=\spann\set{X_1, X_2, X_3}\end{equation*} are rational, then there
exists a strong Malcev basis $\bbase'$ of $\mathfrak g$ strongly based
on $\Gamma$ and passing through $\mathfrak z(\mathfrak g)$ and $\derived$.
Let
\begin{equation*}
\mbox{\sans{P}}_{_{\!\!\!\bbase \rightarrow \bbase'}} =\left(
                                           \begin{array}{ccccc}
                                             a_{11} & a_{12} & a_{13} & a_{14} & a_{15} \\
                                             a_{21} & a_{22} & a_{23} & a_{24} & a_{25} \\
                                             0 & 0 & a_{33} & a_{34} & a_{35} \\
                                             0 & 0 & 0 & a_{44} & a_{45}\\
                                             0 & 0 & 0 & a_{54}& a_{55} \\
                                           \end{array}
                                         \right)
\end{equation*}
be the change of basis matrix from the basis $\bbase$ to the basis
$\bbase'$. Let $\Phi_*\in \mathrm{{Aut}}(\mathfrak g_{5,4})$ def\/ined by
\begin{gather*}
  \Phi_*(X_5)  =  a_{55} X_5+a_{45} X_4+a_{35} X_3+a_{25} X_2+a_{15} X_1, \\
  \Phi_*(X_4)  =  a_{54} X_5+a_{44}X_4+ a_{34} X_3+a_{24} X_2+a_{14}
X_1.
\end{gather*}
Then it is  not hard to verify that there exist $b_{11}, b_{21},
b_{12}, b_{22}, b_{13}, b_{23},
 b_{33}\in \dsr$ such that
\begin{equation*} \Gamma \simeq \Phi^{-1}(\Gamma)
= \exp(\dsz e_1)\exp(\dsz e_2) \exp(\dsz e_3) \exp(\dsz X_4)
\exp(\dsz X_5),
\end{equation*} where $e_1=b_{11}X_1+b_{21}X_2$,
$e_2=b_{12} X_1+b_{22}X_2$, $e_3 =b_{13} X_1+b_{23}X_2+b_{33}X_3$. On
the other hand, as $\Phi^{-1}(\Gamma)$ is a subgroup of $G$, then we
have \begin{equation*}%\label{33}
\exp(X_5) \exp(X_4)
\exp(-X_5)\in\Phi^{-1}(\Gamma).\end{equation*} Then there exist
integer coef\/f\/icients $t_1$, $t_2$, $t_3$, such that
\[
\exp(X_5)
\exp(X_4) \exp(-X_5)=\exp(t_1 e_1) \exp(t_2e_2)
\exp(t_3e_3)\exp(X_4).
\] The above equation may obviously
rewritten as
\begin{equation}\label{om1}
X_3+\frac{1}{2}X_2 =t_1e_1+t_2e_2+t_3e_3-\frac{1}{2}t_3b_{33}X_1.
\end{equation}
Using a similar technique, we obtain that there exist integer
coef\/f\/icients $x_1$, $x_2$, $y_1$, $y_2$ satisfy $x_1^2+x_2^2\ne 0$ and
$y_1^2+y_2^2\ne 0$ such that
\begin{equation}\label{om2}
b_{33}X_2=x_1e_1+x_2e_2
\end{equation}
and
\begin{equation}\label{om3}
b_{33}X_1=y_1e_1+y_2e_2.
\end{equation}
From the equations (\ref{om1}), (\ref{om2}) and (\ref{om3}), it is
clear that the coef\/f\/icients $b_{ij}$ belong to $\dsq$. Let $m$ be
the least common multiple of the denominators of the rational
numbers $b_{ij}$ and let $\pi_*\in \mathrm{{Aut}}(\mathfrak g_{5,4})$
such that $\pi_*(X_5)= X_5$ and $\pi_*(X_4)= m X_4$. Then, we obtain
\begin{gather*}
\Gamma\simeq \exp(\dsz(m^2b_{11}X_1+m b_{21}X_2))
\exp(\dsz(m^2b_{12}X_1+m b_{22}X_2))\\
\phantom{\Gamma\simeq}{} \times\exp(\dsz(m^2b_{13}X_1+m
b_{23}X_2+mb_{33}X_3))\exp(\dsz(mX_4))\exp(\dsz X_5).
\end{gather*}
By Theorem \ref{Hermite}, let
\[
D=\left(
  \begin{array}{ccc}
    c_{11} & c_{12} & c_{13} \\
    0 & c_{22} & c_{23} \\
    0 & 0 & c_{33} \\
  \end{array}
\right)
\]
 be the Hermite normal form of
\[
\left(
    \begin{array}{ccc}
      m^2b_{11} & m^2b_{12} & m^2b_{13} \\
      m b_{21} & m b_{22} & m
b_{23} \\
      0 & 0 & mb_{33} \\
    \end{array}
  \right).
\]
It follows that
\begin{gather*}
\Gamma\simeq \Gamma_{\llbracket D, m  \rrbracket}=
\exp(\dsz(c_{11}X_1)) \exp(\dsz(c_{12}X_1+c_{22}X_2))\\
\phantom{\Gamma\simeq \Gamma_{\llbracket D, m  \rrbracket}=}{}\times \exp(\dsz(c_{13}X_1+c_{23}X_2+c_{33}X_3))\exp(\dsz(mX_4))\exp(\dsz
X_5).
\end{gather*}
It remains to prove that the block matrix  $\llbracket D, m
\rrbracket = \left(
               \begin{array}{cc}
                 D & 0 \\
                 0 & m \\
               \end{array}
             \right)
\in \enss_4$. We need only to prove~(\ref{Mariem1}) and~(\ref{Mariem2}). Since
\begin{equation*}%\label{Khira}
\exp(X_5)\exp(mX_4)\exp(-X_5)\in \Gamma_{\llbracket D, m
\rrbracket}\end{equation*}
 then
 there exist integer coef\/f\/icients $a_{14}$, $a_{24}$, $a_{34}$, such that
\begin{gather}
e^{\ad{X_5}}(mX_4)= a_{14} (c_{11}X_1)
+a_{24}(c_{12}X_1+c_{22}X_2)\nonumber\\
\phantom{e^{\ad{X_5}}(mX_4)=}{} +
a_{34}(c_{13}X_1+c_{23}X_2+c_{33}X_3)+ mX_4-\frac{1}{2}a_{34}
c_{33}m X_1.\label{c4}
\end{gather} Similarly, we establish
\begin{gather}
e^{\ad{X_5}}(c_{13}X_1+c_{23}X_2+c_{33}X_3)=
a_{13} (c_{11}X_1)+a_{23}(c_{12}X_1+c_{22}X_2)\nonumber\\
\phantom{e^{\ad{X_5}}(c_{13}X_1+c_{23}X_2+c_{33}X_3)=}{}+c_{13}X_1+c_{23}X_2+c_{33}X_3,\label{c3}
\\
\label{c2}e^{\ad{X_5}}(c_{12}X_1+c_{22}X_2)=
a_{12} (c_{11}X_1)+(c_{12}X_1+c_{22}X_2),\\
\label{c1}e^{\ad{X_5}}(c_{11}X_1)=
c_{11}X_1.
\end{gather} We see that the conditions (\ref{c4})--(\ref{c1}) boil down to the matrix
equation
\begin{gather*}
\left(\begin{array}{cccc}
1 & 0&0&0\\
0 & 1&1&\frac{1}{2} \\
0&0&1&1\\
0&0&0&1
\end{array} \right)\llbracket D, m  \rrbracket =  \left(
                            \begin{array}{cccc}
                              1 & 0 & -\frac{m}{2} & 0 \\
                              0 & 1 & 0 & 0 \\
                              0 & 0 & 1 & 0 \\
                              0 & 0 & 0 & 1
                            \end{array}
                          \right) \llbracket D, m  \rrbracket \left(\begin{array}{cccc}
                              1 & a_{12} & a_{13} & a_{14} \\
                              0 & 1 & a_{23} & a_{24} \\
                              0 & 0 & 1 & a_{34} \\
                              0 & 0 & 0 & 1
                            \end{array}
                          \right).
\end{gather*}
This proves (\ref{Mariem1}). The proof of the second equality
(\ref{Mariem2}) is shown similarly (repeat the proof of
(\ref{Mariem1}) verbatim, using $m X_4$ instead of $X_5$). The proof
of theorem is complete.
\end{proof}

\section[Uniform subgroups of  $G_{5, 6}$]{Uniform subgroups of  $\boldsymbol{G_{5, 6}}$}\label{section7}

Let $\mathfrak g_{5,6}$ be the three-step nilpotent Lie algebra with
basis
\[
\bbase=(X_1, \ldots, X_5)
\] with Lie brackets are given by
\begin{equation}\label{crochet6-6}
[X_5, X_4]=X_3,\qquad [X_5, X_3]=X_2,\qquad [X_5, X_2]=X_1,\qquad [X_4,
X_3]=X_1
\end{equation} and the non-def\/ined brackets being equal to
zero or obtained by antisymmetry. Let $G_{5, 6}$ be the simply
connected Lie group with Lie algebra $\mathfrak g_{5,6}$.
 Let $\ens_6$ be the set of all invertible integer matrices
of the form
\[
\left(
    \begin{array}{cccc}
      \alpha_{11} & \alpha_{12} & \alpha_{13} & 0 \\
      0 & \alpha_{22} & \alpha_{23} & 0 \\
      0 & 0 & \alpha_{33} & 0 \\
      0 & 0 & 0 & \alpha_{44} \\
    \end{array}
  \right).
  \]
We denote by $\enss_6$ the subset of $\mathrm{Mat}(5, \dsz)$
consisting of all invertible integer matrices of the form
\[
\llbracket D, m \rrbracket = \left(
               \begin{array}{cc}
                 D & 0 \\
                 0 & m \\
               \end{array}
             \right)
             \]
satisfying
\begin{equation*}%\label{Imrano}
 D^{-1} \left(
           \begin{array}{cccc}
             \dfrac{\alpha_{44} m^3}{6}+ \dfrac{\alpha_{44} ^2}{2} & m \alpha_{23}
+\dfrac{m^2\alpha_{33}}{2} & m\alpha_{22} & \alpha_{44} \alpha_{23} \\
             \dfrac{\alpha_{44} m^2}{2} & m \alpha_{33} & 0 & 0 \\
             \alpha_{44} & 0 & 0 & 0 \\
             0 & 0 & 0 & 0 \\
           \end{array}
         \right)\in \mathrm{Mat}(4, \dsz),
\end{equation*}
where $m\in \dsn^*$, the block matrix $D=(\alpha_{ij}:\, 1\leq i,
j\leq 4)\in \ens_6$. Let $\sim$ be the equivalence relation on
$\enss_6$ given by
\begin{gather*}
\llbracket D, m \rrbracket \sim \llbracket D', m' \rrbracket \iff
\llbracket D, m \rrbracket \llbracket D', m' \rrbracket^{-1} =
\mathrm{diag}[a^5, a^4, a^3, a^2,a],\qquad \mbox{for some} \quad a\in
\dsq^*.
\end{gather*}

  We now state the f\/inal result of this paper.
\begin{theorem}\label{uniform(5,6)} With the above notation, we have
\begin{itemize}\itemsep=0pt
\item[$1.$]  If $\llbracket D, m  \rrbracket\in \enss_6$, then
\begin{equation*}\label{uniform-g5-2}
\Gamma_{\llbracket D, m  \rrbracket} = \exp( \dsz e_1) \exp(\dsz
e_2)  \exp(\dsz e_3) \exp(\dsz e_4) \exp(\dsz e_5),
\end{equation*}
 where the vectors $e_j$ $(1\leq j\leq 5)$ are the column vectors of
 $\llbracket D, m  \rrbracket$ in the basis $(X_1, \ldots, X_5)$,
 is a  uniform
subgroup of $G_{5,6}$.
 \item[$2.$] If $ \Gamma$ is a uniform subgroup of $G_{5,6}$, then there exist
 $\llbracket D, m  \rrbracket\in \enss_6$
and $\Phi \in \mathrm{{Aut}}(G_{5,6})$ such that
  $\Phi(\Gamma)= \Gamma_{\llbracket D, m  \rrbracket}$.
\item[$3.$] For $\llbracket D, m  \rrbracket, \llbracket D', m'  \rrbracket\in
  \enss_6$, the subgroups $\Gamma_{\llbracket D, m  \rrbracket}$ and
$\Gamma_{\llbracket D', m'  \rrbracket}$ are isomorphic if and only
if $\llbracket D, m  \rrbracket\sim\llbracket D', m'  \rrbracket$.
 \end{itemize}
\end{theorem}
  We f\/irst state and prove two lemmas.
\begin{lemma} The following ideals
\[
\mathfrak a_i=\spann\set{X_1, \ldots, X_i}\qquad (1\leq i\leq 4)
\]
are rational  with respect to any rational structure on $\mathfrak
g_{5,6}$.
\end{lemma}
\begin{proof} It is clear that $\mathfrak a_1= \mathfrak z(\mathfrak g_{5, 6})$, $\mathfrak
a_2=\mathscr{C}_2(\mathfrak g_{5, 6})$ and $\mathfrak a_3=\EuScript{D}(\mathfrak
g_{5, 6})$. Then the rationality of~$\mathfrak a_1$, $\mathfrak a_2$, $\mathfrak
a_3$ follows from Proposition~\ref{rational-ascending-central-series} and Proposition~\ref{rational-descending-central-series}. On the other hand, we have
\[
\cent(\mathfrak a_2)= \mathfrak a_4.
\]
We conclude from  Proposition \ref{rat-centralizer}  that $\mathfrak
a_4$ is rational.
\end{proof}

\begin{lemma}\label{00} We have
\begin{gather*}
\mathrm{{Aut}}(\mathfrak g_{5,6}) = \Bigg\{\Phi\in \mathrm{{End}}(\mathfrak
g_{5,6}):\ \mathrm{{Mat}}(\Phi, \bbase) =\left(
                                           \begin{array}{ccccc}
                                             a_{55}^5 & a_{12} & a_{13} & a_{14} & a_{15} \\
                                             0 & a_{55}^4 & a_{23} & a_{24} & a_{25} \\
                                             0 & 0 & a_{55}^3 & a_{34} & a_{35} \\
                                            0 & 0 & 0 & a_{55}^2 & a_{45}\\
                                             0 & 0 & 0 & 0 & a_{55} \\
                                           \end{array}
                                         \right)\in \mathrm{GL}(5,
\dsr),\\
\hphantom{\mathrm{{Aut}}(\mathfrak g_{5,6}) = \Bigg\{}{} a_{23} = a_{55} a_{34},\  a_{13}= a_{55} a_{24} + a_{45}
a_{34}-a_{35} a_{55}^2,\  a_{12} = a_{45} a_{55}^3 + a_{34}
a_{55}^2\Bigg\}.
\end{gather*}
\end{lemma}
\begin{proof} Let $\Phi \in \mathrm{{Aut}}(\mathfrak g_{5,6})$.
 Since $\mathfrak a_2$ is invariant under~$\Phi$,  then also $\cent(\mathfrak
a_2)=\mathfrak a_4 $ is invariant under~$\Phi$. It follows that the
matrix $\mathrm{{Mat}}(\Phi, \bbase)$ of $\Phi$ has the following
form
\[
\mathrm{{Mat}}(\Phi, \bbase) =\left(
                                           \begin{array}{ccccc}
                                             a_{11} & a_{12} & a_{13} & a_{14} & a_{15} \\
                                             0 & a_{22} & a_{23} & a_{24} & a_{25} \\
                                             0 & 0 & a_{33} & a_{34} & a_{35} \\
                                             0 & 0 & 0 & a_{44} & a_{45}\\
                                             0 & 0 & 0 & 0 & a_{55}
                                           \end{array}
                                         \right)\in \mathrm{GL}(5,
\dsr).
\]
 The remainder of the proof follows from (\ref{crochet6-6}).
\end{proof}
\begin{proof}[Proof of Theorem \ref{uniform(5,6)}] Assertion $1.$ is obvious.
 To prove the second, let $\Gamma$ be a uniform subgroup of
 $G_{5,6}$.
 Since for every $i=1, \ldots, 4$, the ideal $\mathfrak
a_i=\spann\set{X_1, \ldots, X_i}$ is rational, then there exists a
strong Malcev basis $\bbase'$ of $\mathfrak g$ strongly based on
$\Gamma$ and passing through $\mathfrak a_1, \ldots, \mathfrak a_3$ and
$\mathfrak a_4$. Let
\begin{equation*}
\mbox{\sans{P}}_{_{\!\!\!\bbase \rightarrow \bbase'}} =\left(
                                           \begin{array}{ccccc}
                                             a_{11} & a_{12} & a_{13} & a_{14} & a_{15} \\
                                             0 & a_{22} & a_{23} & a_{24} & a_{25} \\
                                             0 & 0 & a_{33} & a_{34} & a_{35} \\
                                             0 & 0 & 0 & a_{44} & a_{45}\\
                                             0 & 0 & 0 & 0 & a_{55} \\
                                           \end{array}
                                         \right)
\end{equation*}
be the change of basis matrix from the basis $\bbase$ to the basis
$\bbase'$. Let $\Phi_*\in \mathrm{{Aut}}(\mathfrak g_{5,6})$ def\/ined by
\begin{gather*}
  \Phi_*(X_5)  =  a_{55} X_5+a_{45} X_4+\cdots+a_{15} X_1, \\
  \Phi_*(X_4)  =  a_{55}^2 X_4+ \frac{a_{55}^2}{a_{44}}(a_{34} X_3+a_{24} X_2+a_{14}
X_1).
\end{gather*}
Then there exist $b_{11}, b_{12}, b_{22}, b_{13}, b_{23},
 b_{33},
  b_{44}\in \dsr$ such that
\begin{equation*}
\Gamma \simeq \Phi^{-1}(\Gamma)
= \exp(\dsz \epsilon_1)\exp(\dsz \epsilon_2) \exp(\dsz \epsilon_3)
\exp(\dsz \epsilon_4) \exp(\dsz X_5),
\end{equation*} where
$\epsilon_1=b_{11}X_1$, $\epsilon_2=b_{12} X_1+b_{22}X_2$, $\epsilon_3
=b_{13} X_1+b_{23}X_2+b_{33}X_3$, $\epsilon_4=b_{44} X_4$. By similar
arguments to those used in (\ref{c4})--(\ref{c1}) we derive the
following system
\begin{gather*}
  b_{44}  =  x_1 b_{33}, \\
  \frac{1}{2} b_{44}  =  x_1 b_{23}+ x_2 b_{22},\\
  \frac{1}{6} b_{44}  =  x_1 b_{13}+ x_2 b_{12}+x_3 b_{11}-\frac{1}{2} x_1b_{44}b_{33},\\
  b_{33}  =  y_1 b_{22}, \\
  b_{23}+\frac{1}{2} b_{33}  =  y_1 b_{12}+ y_2 b_{11},\\
  b_{44} b_{33}   =  z b_{11}, \\
  b_{22}  =  t b_{11},
\end{gather*}
where $x_1,  y_1,  z, t \in \dsn^*$, $x_2, x_3,y_2 \in \dsz$.  By
calculation one has
\[
( b_{11}, b_{12}, b_{22}, b_{13}, b_{23},
 b_{33},
  b_{44},1)\in \mathscr{G}.
\]
  where
\begin{gather*}
\mathscr{G}  =\Bigg\{(\alpha, m)\in \dsq^7\times \dsn^*:\ \alpha=
(\alpha_{1},
\alpha_{2},\alpha_{3},\alpha_{4},\alpha_{5},\alpha_{6},\alpha_{7})\
\mbox{ such that }
\\
\phantom{\mathscr{G}  =\Bigg\{}{}
\alpha_{1}= \frac{b m^4}{zs^2a^2}, \ \alpha_{2}= \frac{b m^3}{2zs^2a}-\frac{yb m^4}{z^2s^3a}+\frac{b
m^4}{2zs^2a} -\frac{t b m^4}{zs^3a^2}, \ \alpha_{3}= \frac{b m^3}{zs^2a},\\
\phantom{\mathscr{G}  =\Bigg\{}{} \alpha_{4}= \frac{b m^2}{6zsa}-\frac{y b m^3}{2z^2s^2a}+\frac{y^2b
m^4}{z^3s^3a}-\frac{y b m^4}{2z^2s^2a}
+\frac{y t b m^4}{z^2s^3a^2}-\frac{x b m^4}{z^2s^2a^2}+\frac{b^2 m^4}{2zs^2a^2},\\
\phantom{\mathscr{G}  =\Bigg\{}{} \alpha_{5}= \frac{b m^2}{2zsa}-\frac{yb m^3}{z^2s^2a}, \
 \alpha_{6}= \frac{b m^2}{zsa}, \ \alpha_{7}= \frac{b m^2}{sa},\  a, b, s, z \in \dsz^*,\ x, y, t
\in \dsz\Bigg\}.
\end{gather*}
From this we deduce that the coef\/f\/icients $b_{ij}$ belong to $\dsq$.
Let $m$ be the least common multiple of the denominators of the
rational numbers $b_{ij}$. Let $\pi\in \mathrm{Aut}(G_{5,6})$ such
that
\[
\mathrm{Mat}(\pi_*, \bbase) = \texttt{diag}\,[m^5, m^4, m^3, m^2,
m].
\] Then
\begin{equation*} \Gamma \simeq \pi(\Phi^{-1}(\Gamma))
= \exp(\dsz e_1)\exp(\dsz e_2) \exp(\dsz e_3) \exp(\dsz e_4)
\exp(m\dsz X_5)\end{equation*} such that the matrix $[e_1, \ldots,
e_4]$ with column vectors $e_1, \!\ldots, \!e_4$ expressed in the basis
$\{X_1, \!\ldots, \!X_4\}$, belongs to $\ens_6$. Write
\[
[e_1, \ldots,
e_4]= \left(
    \begin{array}{cccc}
      \alpha_{11} & \alpha_{12} & \alpha_{13} & 0 \\
      0 & \alpha_{22} & \alpha_{23} & 0 \\
      0 & 0 & \alpha_{33} & 0 \\
      0 & 0 & 0 & \alpha_{44}
    \end{array}
  \right).
\]
  By  similar techniques as above we can prove that
$(\alpha_{11}, \alpha_{12}, \alpha_{13}, \alpha_{22}, \alpha_{23},
\alpha_{33},\alpha_{44}, m)\in \mathscr{G}$. This equivalent that
\[
[e_1, \ldots,
e_4]^{-1}\left(
           \begin{array}{cccc}
             \dfrac{\alpha_{44} m^3}{6}+ \dfrac{\alpha_{44} ^2}{2} & m \alpha_{23}
+\dfrac{m^2\alpha_{33}}{2} & m\alpha_{22} & \alpha_{44} \alpha_{23} \\
             \dfrac{\alpha_{44} m^2}{2} & m \alpha_{33} & 0 & 0 \\
             \alpha_{44} & 0 & 0 & 0 \\
             0 & 0 & 0 & 0
           \end{array}
         \right)\in \mathrm{Mat}(4, \dsz).
\]
 Finally, we achieve with the proof of $3$. We show both directions. Let $a\in \dsq^*$
 such that
\[
\llbracket D, m \rrbracket \llbracket D', m' \rrbracket^{-1} =
\mathrm{diag}\, [a^5, a^4, a^3, a^2,a].
\]
 Consider the linear mapping
$\phi_*: \mathfrak g_{5,6}\longrightarrow \mathfrak g_{5,6}$ def\/ined by
\[
\mathrm{Mat}(\phi_*, \bbase) = \mathrm{diag}\, [a^5, a^4, a^3,
a^2,a].
\]
 It is clear that $\phi\in \mathrm{Aut}(G_{5,6})$ (see
Lemma \ref{00}) and $\phi(\Gamma_{\llbracket D, m \rrbracket})=
\Gamma_{\llbracket D', m' \rrbracket}$. Conversely, suppose there
exists $\phi: \Gamma_{\llbracket D, m \rrbracket} \longrightarrow
\Gamma_{\llbracket D', m' \rrbracket}$  an isomorphism between
$\Gamma_{\llbracket D, m \rrbracket}$ and $\Gamma_{\llbracket D', m'
\rrbracket}$. By Theorem~\ref{extended}, $\phi$ has an extension
$\overline{\phi}\in \Aut(G_{5,6})$. As $\overline{\phi}_*(m X_5)= m'
X_5$ and $\overline{\phi}_*(\alpha_{44} X_4)= \alpha'_{44} X_4$
where the $\alpha_{ij}$ (resp. $\alpha'_{ij}$) are the entries of
$D$ (resp. of $D'$), then by Lemma \ref{00}, the matrix of
$\overline{\phi}_*$ in the basis $\bbase$ has the following form
\[
\mathrm{Mat}(\overline{\phi}_*, \bbase) = \mathrm{diag}\, \big[a^5, a^4, a^3,
a^2,a\big]
\] for some $a\in \dsq^*$. Consequently, we obtain
\[
\llbracket D, m \rrbracket =
\mathrm{diag}\,\big[a^5, a^4, a^3, a^2,a\big]\, \llbracket D', m' \rrbracket.
\]
 This completes the proof.
\end{proof}

\begin{remark} In the statement 2 of Theorem \ref{uniform(5,6)},
 the element $\llbracket D, m \rrbracket$ of $\enss_6$ is not unique.\end{remark}

\subsection*{Acknowledgements}
The authors would like to thank the anonymous
  referees for their critical and valuable comments.

\pdfbookmark[1]{References}{ref}
\LastPageEnding

\end{document}